\documentclass{amsart}
\usepackage{latexsym, amssymb}
\newtheorem{thm}{Theorem}

\newtheorem{cor}{Corollary}

\newtheorem{lem}{Lemma}
\newtheorem{defn}{Definition}

\begin{document}
\begin{center}
{\bf FINITE TYPE LINK HOMOTOPY INVARIANTS}\\
\vspace{.2in}
{\footnotesize BLAKE MELLOR}\\
{\footnotesize Honors College}\\
{\footnotesize Florida Atlantic University}\\
{\footnotesize 5353 Parkside Drive}\\
{\footnotesize Jupiter, FL  33458}\\
{\footnotesize\it  bmellor@fau.edu}\\
 
\vspace{1in}
{\footnotesize ABSTRACT}\\
{\ }\\
\parbox{4.5in}{\footnotesize \ \ \ \ \ In \cite{bn2}, Bar-Natan used unitrivalent diagrams
to show that finite type invariants classify string links up to homotopy.
In this paper, I will construct the correct spaces of chord diagrams and unitrivalent diagrams
for links up to homotopy.  I will use these spaces to show that, far from classifying links
up to homotopy, the only rational finite type invariants of link homotopy are the linking
numbers of the components.
}\\
\end{center}
\input{vpsfig.sty}
\tableofcontents
\section{Introduction} \label{S:intro}
We will begin with a brief overview of finite type invariants.
In 1990, V.A. Vassiliev introduced the idea of {\it Vassiliev} or {\it finite type} knot
invariants, by looking at certain groups associated with the cohomology of the space of
knots.  Shortly thereafter, Birman and Lin~\cite{bl} gave a combinatorial description
of finite type invariants.  We will give a brief overview of this combinatorial theory.
For more details, see Bar-Natan~\cite{bn1}.
\subsection{Singular Knots and Chord Diagrams} \label{SS:chord}
Recall that, in the most general sense, a {\it knot invariant} is a map from
the set of equivalence classes of knots under isotopy to another set
$G$.  We will need to have some additional structure on $G$.  For our purposes,
$G$ will always be at least an associative, commutative ring with an identity.
We first note that we can extend any knot invariant to an invariant of
{\it singular} knots, where a singular knot is an immersion of $S^1$ in
3-space which is an embedding except for a finite number of isolated double
points.  Given a knot invariant $v$, we extend it via the relation:
$$\psfig{file=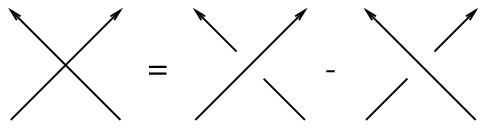}$$
An invariant $v$ of singular knots is then said to be of {\it finite type}, if
there is an integer $d$ such that $v$ is zero on any knot with more than $d$
double points.  $v$ is then said to be of {\it type} $d$.  We denote by $V_d(G)$
the $G$-module generated by $G$-valued finite type invariants of type $d$.  We
can completely understand the space of $G$-valued finite type invariants by understanding
all of the $G$-modules $V_d(G)/V_{d-1}(G)$.  An element of this module is
completely determined by its behavior on knots with exactly $d$ singular
points.  Since such an element is zero on knots with more than
$d$ singular points, any other (non-singular) crossing of the knot can be
changed without affecting the value of the invariant.  This means that
elements of $V_d/V_{d-1}$ can be viewed as functions on the space of
{\it chord diagrams}:
\begin{defn}
A {\bf chord diagram of degree d} is an oriented circle, together with $d$
chords of the circles, such that all of the $2d$ endpoints of the chords are
distinct.  The circle represents a knot, the endpoints of a chord represent
2 points identified by the immersion of this knot into 3-space.
\end{defn}
Functions on the space of chord diagrams which are derived from knot
invariants will satisfy certain relations.  This leads us to the definition
of a {\it weight system}:
\begin{defn}
A $G$-valued {\bf weight system of degree d} is a $G$-valued function $W$ on the
space of chord diagrams of degree $d$ which satisfies 2 relations:
\begin{itemize}
    \item (1-term relation) $$\psfig{file=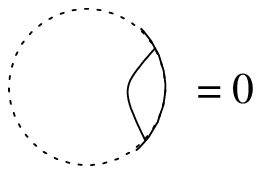}$$
    \item (4-term relation) $$\psfig{file=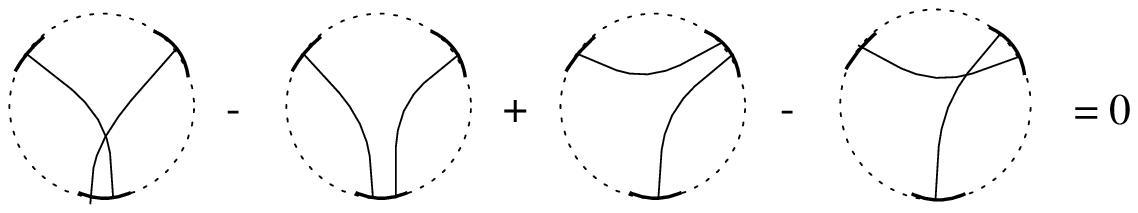}$$
    Outside of the solid arcs on the circle, the diagrams can be anything,
as long as it is the same for all four diagrams.
\end{itemize}
We let $W_d(G)$ denote the space of $G$-valued weight systems of degree d.
\end{defn}
From now on, we will be considering the case when $G = {\mathbb R}$.  We will
simplify our notation by letting $V_d = V_d({\mathbb R})$ and $W_d = W_d({\mathbb R})$.
Bar-Natan~\cite{bn1} defines maps $w_d: V_d \rightarrow W_d$ and $v_d: W_d \rightarrow V_d$.
$w_d$ is defined by {\it embedding} a chord
diagram $D$ in ${\bf R}^3$ as a singular knot $K_D$, with the chords corresponding to
singularities of the embedding (so there are $d$ singularities).  Any two
such embeddings will differ by crossing changes, but these changes will not effect the value of
a type $d$ Vassiliev invariant on the singular knot.  Then, for any $\gamma \in V_d$,
we define $w_d(\gamma)(D) = \gamma(K_D)$.  Bar-Natan shows that this is, in fact, a weight
system.  The 1-term relation is satisfied because of the first Reidemeister move, and the 4-term
relation is essentially the result of rotating a third strand a full turn around a double point.
Note that this argument will work for any ring $G$, not just $G = {\mathbb R}$.
$v_d$ is much more complicated to define, using the Kontsevich integral.  This does require that
we consider real-valued invariants; it is still an open question whether an analogous prcedure
can be found for other rings (or even finite fields).  For a full treatment
of the Kontsevich integral, see Bar-Natan~\cite{bn1} and Le and Murakami~\cite{lm}.  We will
simply mention the few
facts and properties we need, primarily following Le and Murakami.  Using a Morse function, any knot
(or link or string link) can be decomposed into elementary {\it tangles}:
$$\psfig{file=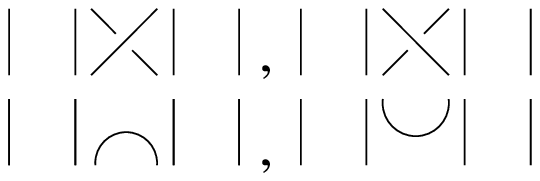}$$
Le and Murakami define a map $Z$ from an elementary tangle with $k$ strands
to the space of chord diagrams on $k$ strands.  This map respects composition of tangles:
if $T_1\cdot T_2$ is the tangle obtained by placing $T_1$ on top of $T_2$, then
$Z(T_1\cdot T_2) = Z(T_1)Z(T_2)$.  Le and Murakami prove that this map gives a real-valued
isotopy invariant of knots and links.
Given a degree $d$ weight system $W$, and a knot $K$, we now define $v_d(W)(K) = W(Z(K))$.
Bar Natan~\cite{bn1} shows that $w_d$ and $v_d$ are ``almost'' inverses.  More precisely,
$w_d(v_d(W)) = W$ and $v_d(w_d(\gamma))-\gamma \in V_{d-1}$.  As a result,
(see \cite{bl,bn1,st}) the space $W_d$ of weight
systems of degree $d$ is isomorphic to $V_d/V_{d-1}$.  For convenience,
we will usually take the dual approach, and simply study the real vector space of
chord diagrams of degree $d$ modulo the 1-term and 4-term relations.
The dimensions of these spaces have been computed for $d \leq
12$ (see Bar-Natan~\cite{bn1} and Kneissler~\cite{kn}).  It is useful to combine all of
these spaces into a graded module via direct sum.  We can give this module a Hopf algebra
structure by defining an appropriate product and co-product:
\begin{itemize}
    \item  We define the product $D_1 \cdot D_2$ of two chord diagrams
$D_1$ and $D_2$ as their connect sum.  This is well-defined modulo the 4-term
relation (see \cite{bn1}).
$$\psfig{file=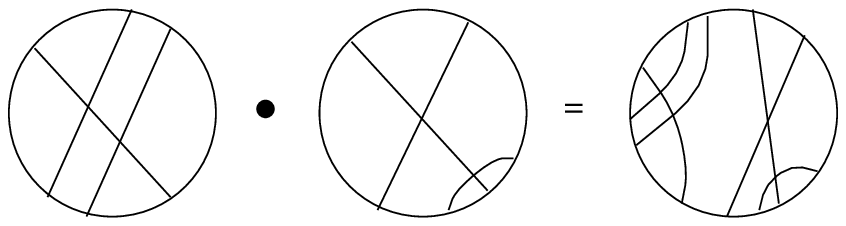}$$
    \item  We define the co-product $\Delta(D)$ of a chord diagram $D$ as
follows:
$$
\Delta(D) = {\sum_J D_J' \otimes D_J''}
$$
where $J$ is a subset of the set of chords of $D$, $D_J'$ is $D$ with all the
chords in $J$ removed, and $D_J''$ is $D$ with all the chords not in J
removed.
$$\psfig{file=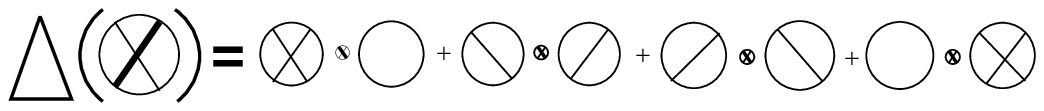}$$
\end{itemize}
It is easy to check the compatibility condition $\Delta(D_1\cdot D_2) = \Delta
(D_1)\cdot\Delta(D_2)$.
\subsection{Unitrivalent Diagrams} \label{SS:ud}
\begin{figure}
    $$\psfig{file=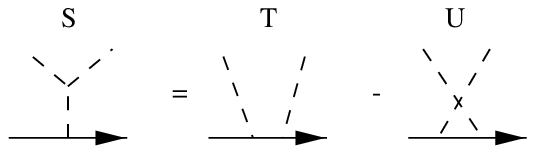}$$
    \caption{STU relation} \label{F:stu}
    \end{figure}

\begin{figure}
    $$\psfig{file=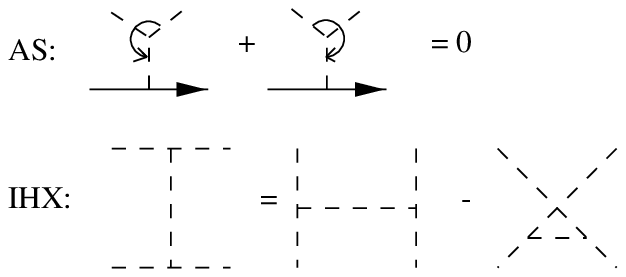}$$
    \caption{AS and IHX relations} \label{F:ihx}
    \end{figure}

It is often useful to consider the Hopf algebra of bounded unitrivalent diagrams, rather
than chord diagrams.  These diagrams, introduced by Bar-Natan~\cite{bn1} (there called
{\it Chinese Character Diagrams}), can
be thought of as a shorthand for writing certain linear combinations of chord
diagrams.  We define a {\it bounded unitrivalent diagram} to be a unitrivalent graph, with
oriented vertices, together with a bounding circle to which all the univalent vertices are
attached.  We also require that each component of the graph have at least one univalent
vertex (so every component is connected to the boundary circle).  We define the space $A$
of bounded unitrivalent diagrams as the quotient of the space of all bounded unitrivalent
graphs by the $STU$ relation, shown in Figure~\ref{F:stu}.
As consequences of $STU$ relation, the anti-symmetry ($AS$) and $IHX$ relations, see
Figure~\ref{F:ihx}, also hold in $A$.
Bar-Natan shows that $A$ is isomorphic to the algebra of chord diagrams.
We can get an algebra $B$ of {\it unitrivalent diagrams} by simply removing the bounding circle
from the diagrams in $A$, leaving graphs with trivalent and univalent vertices, modulo
the $AS$ and $IHX$ relations.  Bar-Natan shows that the spaces $A$ and $B$ are
isomorphic.  The map $\chi$ from $B$ to $A$ takes a diagram to the linear combination of
all ways of attaching the univalent vertices to a bounding circle, divided by total number of
such ways (T. Le noticed that this factor, missing in \cite{bn1}, is necessary to preserve the
comultiplicative structure of the algebras).  The inverse map $\sigma$
turns a diagram into a linear combination of diagrams by performing sequences of ``basic
operations,'' and then removes the bounding circle.  The two basic operations are:
$$\psfig{file=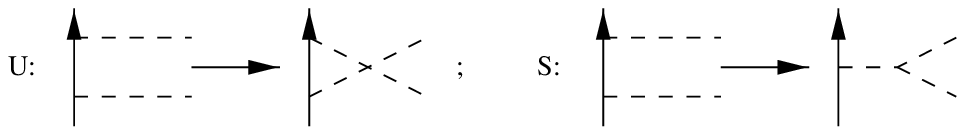}$$
\section{String Links, Links and Homotopy} \label{S:string}
\subsection{String Links} \label{SS:string_links}
Bar-Natan~\cite{bn2} extended the theory of finite type invariants to string links.
\begin{defn} (see Habegger and Lin~\cite{hl})
Let D be the unit disk in the plane and let I = [0,1] be the unit interval.  Choose k points
$p_1,..., p_k$ in the interior of D, aligned in order along the the x-axis.  A {\bf string
link} $\sigma$ of k components is a smooth proper imbedding of k disjoint copies of I into
$D \times I$:
$$\sigma:\ \bigsqcup_{i=1}^k{I_i} \rightarrow D \times I$$
such that $\sigma|_{I_i}(0) = p_i \times 0$ and $\sigma|_{I_i}(1) = p_i \times 1$.  The image
of $I_i$ is called the ith string of the string link $\sigma$.
\end{defn}
Essentially, everything works the same way for string links as for knots.
The bounding circle of the bounded unitrivalent diagrams now becomes a set of
bounding line segments, each labeled with a color, to give an algebra $A^{sl}$ (the multiplication
is given by placing one diagram on top of another).  The univalent diagrams are unchanged,
except that each univalent vertex is also labeled with a color to
give the space $B^{sl}$.  The isomorphisms $\chi$ and $\sigma$ between $A$ and $B$ easily
extend to isomorphisms $\chi^{sl}$ and $\sigma^{sl}$ between $A^{sl}$ and $B^{sl}$, just
working with each color separately.  In addition, there are obvious maps $w_d^{sl}$ and
$v_d^{sl}$ analogous to $w_d$ and $v_d$ (we just need to keep track of colors).
\subsection{Links} \label{SS:links}
\begin{figure}
    $$\psfig{file=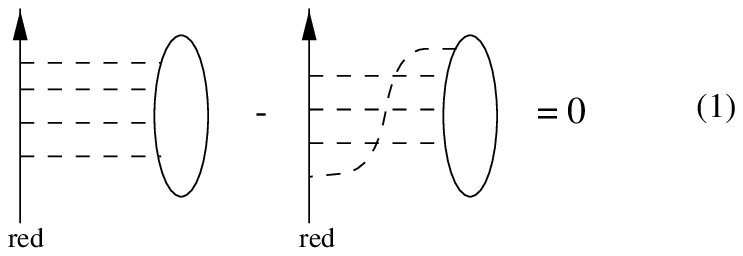}$$
    \caption{The link relation for chord diagrams} \label{F:link_rel}
    \end{figure}

\begin{figure}
    $$\psfig{file=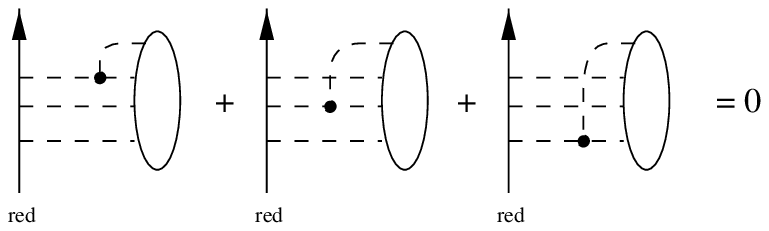}$$
    \caption{The link relation for Bounded Unitrivalent Diagrams} \label{F:link_rel2}
    \end{figure}

\begin{figure}
    $$\psfig{file=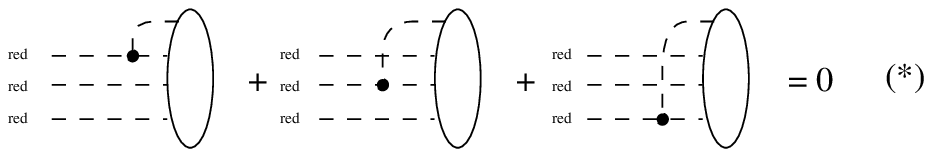}$$
    \caption{The link relation for unitrivalent diagrams} \label{F:link_rel3}
    \end{figure}

The obvious definition of chord diagrams for links is simply to replace the bounding line
segments with bounding circles.  However,
these diagrams are difficult to work with, and it is unclear how to define the
unitrivalent diagrams.  Unlike for a knot, closing up the components of a string link of several
components is not a trivial operation, so we need to impose some relations on the space
of unitrivalent diagrams.
Since we understand the spaces of chord diagrams and unitrivalent diagrams for string links, it
would be useful to be able to express these spaces for links as quotients of the spaces for
string links.  The question is then, what relations do we need?  One relation is fairly
obvious.  When we construct the space $A^l$ of bounded unitrivalent diagrams for links,
we replace the bounding line segments of $A^{sl}$ with directed circles.  Bar-Natan et. al.
observed (see Theorem 3,~\cite{bgrt}) that this is exactly
equivalent to saying that the ``top'' edge incident to one of the line segments
can be brought around the circle to be on the ``bottom.''  So we can write $A^l$ as the
quotient of $A^{sl}$ by relation (1), shown in Figure~\ref{F:link_rel} (where the figure
shows {\it all} the chords with endpoints on the red component).
Then the Kontsevich integral for links, $Z^l$, is defined by cutting the link to make a string
link, applying the Kontsevich integral for string links, and then taking the quotient by
relation (1).  Now $w_d^l$ and $v_d^l$ are defined similarly to $w_d$ and $v_d$.
Given a link invariant $\gamma$ of type $d$ and a diagram $D$ of degree $d$ in $A^l$,
$w_d^l(\gamma)(D) = \gamma(L_{\hat{D}})$, where $\hat{D}$ is the closure of the
diagram $D$ (i.e. the bounding line segments are closed to circles).
$L_{\hat{D}}$ is well-defined by Theorem 3 of
\cite{bgrt}.  Defining $v_d^l$ is even easier, now that we have $Z^l$.  Given a weight system
(element of the graded dual of $A^l$) $W$ and a link $L$, we define
$v_d^l(W)(L) = W(Z^l(L))$.
One advantage of this formulation of $A^l$ is that it enables us to define
the space $B^l$ of unitrivalent diagrams as a quotient of the already known
space $B^{sl}$.  This was done by Bar-Natan et. al.
Using the $STU$ relation, we can rewrite relation (1) as in Figure~\ref{F:link_rel2}.
This suggests how we should define the space $B^l$.  We will take the quotient of $B^{sl}$
by the relations (*) shown in Figure~\ref{F:link_rel3},
where the univalent vertices shown are {\it all} the univalent vertices of a given color.  With
these definitions, Bar-Natan et. al. proved that $A^l$ and $B^l$ are isomorphic:
\begin{thm} \label{T:isotopy-descends} (Theorem 3, \cite{bgrt})
The isomorphism between $A^{sl}$ and $B^{sl}$ descends to an isomorphism
between $A^l$ and $B^l$.
\end{thm}
\subsection{Link Homotopy} \label{SS:homotopy}
The idea of {\it link homotopy} (or just {\it homotopy}) was introduced by Milnor~\cite{mi1}.
Two links are homotopic if one can be transformed into the other through a sequence of ambient
isotopies of $S^3$ and crossing changes of a component with itself (but {\it not} crossing
changes of different components).  The definition for string links is similar.  Habegger
and Lin~\cite{hl} succeeded in classifying string links and links up to homotopy.
We want to extend the results of the last section to string links and links considered up to
homotopy.  For string links, this has already been done by Bar-Natan~\cite{bn2}.  Bar-Natan
describes the algebras $A^{hsl}$ and $B^{hsl}$ of bounded and unbounded unitrivalent diagrams
for string links up to homotopy.  In brief, we take the quotient of $A^{sl}$
(resp. $B^{sl}$) by the space of {\it boring} diagrams.  A diagram is {\it boring}
if it has (1) two univalent vertices on the same component (resp. assigned the same color),
or (2) non-trivial first homology.  In other words, we are left with tree diagrams
with no more than one univalent vertex on each component (resp. of each color).
Bar-Natan then defines $w_d^{hsl}$ and $v_d^{hsl}$ in the usual way, and shows that
they are ``almost'' inverses in the same sense that $w_d$ and $v_d$ are.
All of this extends to links just as it did for isotopy.  We define $A^{hl}$ as the quotient
of $A^{hsl}$ by relation (1), and $B^{hl}$ as the quotient of $B^{hsl}$ by relation (*).  We
then define $Z^{hl}$, $w_d^{hl}$, and $v_d^{hl}$ just as we did for links up to homotopy.
Finally, the arguments of Bar-Natan et. al. carry through to show:
\begin{thm} \label{T:homotopy-descends} (Theorem 3, \cite{bgrt})
The isomorphism between $A^{hsl}$ and $B^{hsl}$ descends to an isomorphism
between $A^{hl}$ and $B^{hl}$.
\end{thm}
{\bf Remark:}  By results of Habegger and Masbaum (see Remark 2.1 of \cite{hm}), $Z^{hl}$ is
the {\it universal} finite type invariant of link homotopy.  By this we mean that it dominates
all other such invariants.
\section{The Size of $B^{hl}$} \label{S:size}
Now that we have properly defined the space $B^{hl}$ of unitrivalent diagrams for link homotopy, we want
to analyze it more closely.  We will consider the case when $B^{hl}$ is a vector space over the reals
(or, more generally, a module over a ring of characteristic 0).
In particular, we would like to know exactly which diagrams of $B^{hsl}$
are in the kernel of the relation (*) (i.e. are 0 modulo (*)).  We will find that the answer is
``almost everything'' - to be precise, any unitrivalent diagram with a component of degree 2 or more.
We will start by proving a couple of base cases, and then prove the rest of the theorem by induction.
Let $B^{hsl}(k)$ denote the space of unitrivalent diagrams for string link homotopy
with $k$ possible colors for the univalent vertices (i.e. we are looking at links with $k$ components).
Consider a diagram $D \in B^{hsl}(k)$.  Recall from the previous sections that each component of $D$
is a tree diagram with at most one endpoint of each color.  Since a unitrivalent tree with $n$
endpoints has $2n-2$ vertices, and hence degree $n-1$, $D$ cannot have any components of degree greater
than $k-1$.
{\bf Notation:}  Before we continue, we will introduce a bit of notation which will be useful in this
section.  Given a unitrivalent diagram $D$, we define $m(D;i,j)$ to be the number of components
of $D$ which are simply line segments with ends colored $i$ and $j$, as shown below:
$$i-----j$$
\subsection{Base cases} \label{SS:base}
\begin{figure}
    $$\psfig{file=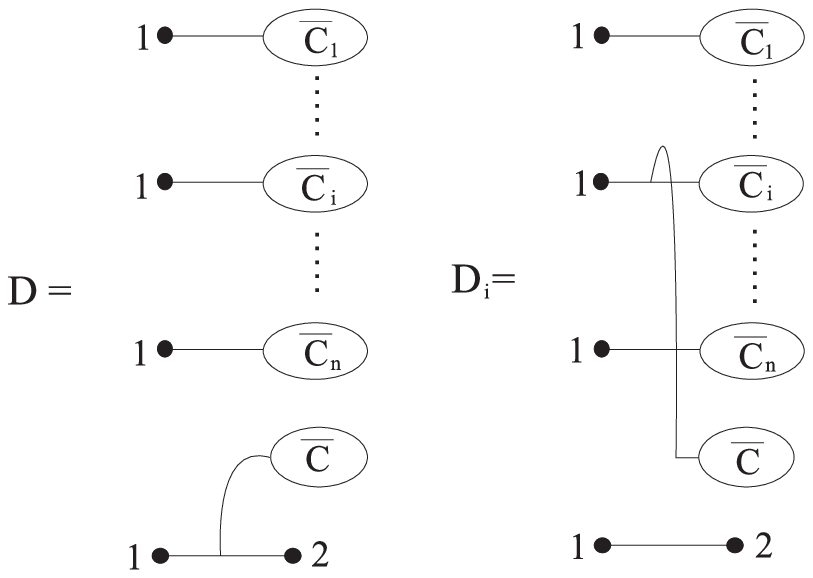}$$
    \caption{Diagrams arising from relation (*)} \label{F:arise}
    \end{figure}

\begin{figure}
    $$\psfig{file=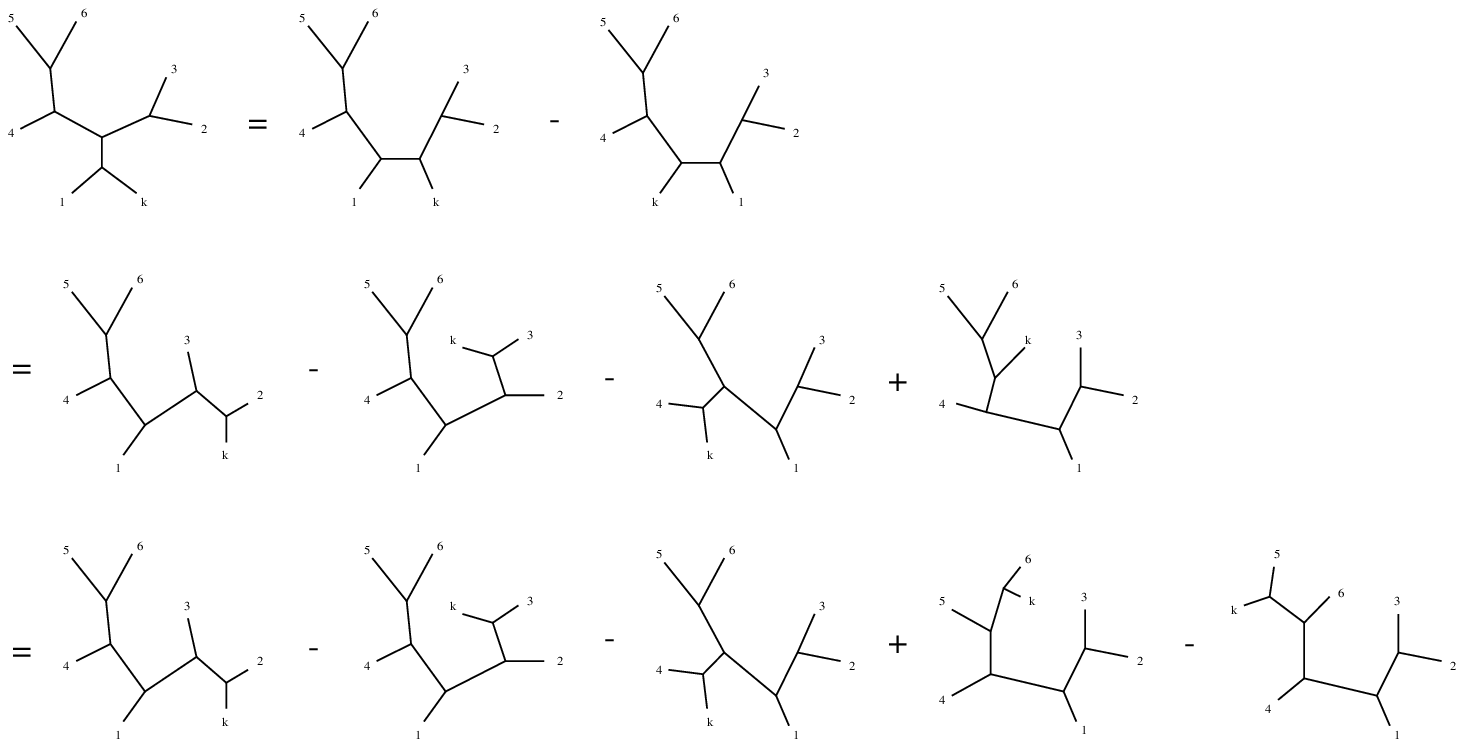}$$
    \caption{Using the IHX relation to decompose a diagram} \label{F:expand}
    \end{figure}

\begin{lem} \label{L:base1}
If D has a component C of degree k-1 (with $k \geq 3$), then D is trivial modulo (*).
\end{lem}
{\sc Proof:}  $C$ has one endpoint of each color $1,2,...,k$.  Without loss of generality, we
may assume that $C$ has a branch as shown, where $\bar{C}$ denotes the remainder of $C$:
$$C:\ \ \begin{matrix} \bar{C} \\ | \\ | \\ 1-----2 \end{matrix}$$
We are going to apply (*) with the color 1.  Let $\{C_1,...,C_n\}$ be the components of $D$ with an
endpoint colored 1.  So, ignoring the other components of $D$, we have the diagrams of
Figure~\ref{F:arise} (where $\bar{C_i}$ denotes all of $C_i$ except for the endpoint colored 1).
(*) then implies that $D + \sum{D_i} = 0$.  If $C_i$ is just a line segment with endpoints colored 1 and
2, then $D_i = D$.  Otherwise, $D_i$ will have a {\it boring} component (since $\bar{C_i}$ will have an
endpoint of some color $j \in {3,...,k}$, and $\bar{C}$ has an endpoint of each color 3,...,$k$, including
$j$, $D_i$ will have a component with two endpoints colored $j$), and hence be trivial in $B^{hsl}$.
Therefore, we find that $D+m(D;1,2)D = 0$ where $m(D;1,2) \geq 0$.
We can divide both sides by $1+m(D;1,2)$ (since we are working over the reals, which have
characteristic 0) to conclude that $D = 0$.  $\Box$
\begin{lem} \label{L:base2}
If D has a component C of degree k-2 (with $k \geq 4$), then D is trivial modulo (*).
\end{lem}
{\sc Proof:}  Without loss of generality, $C$ has endpoints colored $1,2,...,k-1$.  We will
prove the lemma by inducting on $m(D;1,k)$;
inducting among the set of diagrams having a component with endpoints colored $1,2,...,k-1$.
As in the previous lemma, we may assume that $C$ has a branch as shown:
$$C:\ \ \begin{matrix} \bar{C} \\ | \\ | \\ 1-----2 \end{matrix}$$
And conclude that $D+\sum{D_i}=0$, where the $D_i$ are defined as before.  Since $\bar{C}$ contains
endpoints of all colors except 1, 2, and $k$, $D_i$ is boring
unless $C_i$ has one of the following 3 forms (as in Lemma~\ref{L:base1}):
$$(1)\ \ C_i = \ 1-----2$$
$$(2)\ \ C_i = \ 1-----k$$
$$(3)\ \ C_i = \ \begin{matrix} k \\ | \\ | \\ 1-----2 \end{matrix}$$
In the first case, $D_i = D$; and in the second case, $D_i = D'$, where $D'$ is the same as $D$ except
that:
\begin{itemize}
    \item $C$ is replaced by a component $C'$ identical to it except that the endpoint colored
2 in $C$ is colored $k$ in $C'$ (so $\bar{C'} = \bar{C}$).
    \item A line segment with endpoints colored 1 and $k$ has been replaced by a line segment
with endpoints colored 1 and 2.  In other words, $m(D';1,2) = m(D;1,2)+1$ and $m(D';1,k) = m(D;1,k)-1$.
\end{itemize}
In the third case, $D_i$ has a component of
degree $k-1$, and so is 0 modulo (*) by the previous lemma.  Therefore, as in the previous lemma, we find
that $D + m(D;1,2)D + m(D;1,k)D' = 0$.  If $m(D;1,k) = 0$ we conclude, as before, that $D$ is trivial
modulo (*), which proves the base case of our induction.
For the inductive step, we use the IHX relation on $C'$ to decompose $D' = \sum_{i\neq 1,2,k}{\pm D_i'}$,
where $D_i'$ is the same as $D'$ except that $C'$ has been replaced by a component $C_i'$ with endpoints
of the same colors (although arranged differently), and a branch as shown:
$$C_i':\ \ \begin{matrix} \bar{C_i'} \\ | \\ | \\ i-----k \end{matrix}$$
(The decomposition is simply a matter of letting the endpoint colored $k$ ``travel'' the tree -
see Figure~\ref{F:expand} for an example.)  In particular, $m(D_i';a,b) = m(D';a,b)$ for all colors
$a$ and $b$.
We now apply (*) to $D_i'$ using color $i$ (and component $C_i'$), similarly to what we've done
before.  In this case, the only other components which matter (modulo boring diagrams) are
ones which look like one of the following:
$$(1)\ \ i-----k$$
$$(2)\ \ i-----2$$
$$(3)\ \ \begin{matrix} 2 \\ | \\ | \\ i-----k \end{matrix}$$
As before, the first case gives $D_i'$ again, the third case is trivial by Lemma~\ref{L:base1},
and the second case gives a diagram $D_i''$ such that:
\begin{itemize}
    \item $C_i'$ is replaced by a component $C_i''$ identical to it except that the endpoint colored
$k$ in $C_i'$ is colored 2 in $C_i''$ (so $\bar{C_i''} = \bar{C_i'}$).
    \item A line segment with endpoints colored $i$ and 2 has been replaced by a line segment
with endpoints colored $i$ and $k$.  In other words, $m(D_i'';i,k) = m(D_i';i,k)+1$ and
$m(D_i'';i,2) = m(D_i';i,2)-1$.
\end{itemize}
Otherwise, $D_i''$ is the same as $D_i'$; in particular, $m(D_i'';1,k) = m(D_i';1,k) = m(D';1,k) =
m(D;1,k)-1$.
Then (*) tells us that $D_i' + m(D';i,k)D_i' + m(D';2,i)D_i'' = 0$.  Since $D_i''$
has a component of degree $k-2$ with endpoints colored $1,...,k-1$ (namely, $C_i''$), the inductive
hypothesis implies
that $D_i''$ is trivial modulo (*).  Therefore, $(1+m(D';i,k))D_i' = 0$ modulo (*), so $D_i'$ is
trivial modulo (*).  This is true for every $i$, so it immediately follows that $D'$, and hence $D$,
are also trivial modulo (*).  $\Box$
\subsection{Main Theorem} \label{SS:main}
The lemmas of Section~\ref{SS:base} will act as base cases for the main theorem of this section:
\begin{thm} \label{T:induction}
If D has a component C of degree 2 or higher, then D is trivial modulo (*).
\end{thm}
{\sc Proof:}  The method of proof for this theorem is essentially the same as that for
Lemma~\ref{L:base2}.  We will successively apply (*) (and do a single expansion via IHX)
until we obtain a set of diagrams which are all either trivial or repetitions of earlier
diagrams.  We can then backtrack to show that everything disappears.  However, rather than
applying (*) twice, as in Lemma~\ref{L:base2}, we will need to apply it four times.  This
unfortunately makes keeping track of the diagrams somewhat confusing - we have done our best.
Again, the proof is by induction; in this case, it is a nested double induction.  The outer
induction is backwards, on the degree of the largest component of $D$.  The base
cases of this induction are given by Lemma~\ref{L:base1} and Lemma~\ref{L:base2}.  So we assume that
a diagram is trivial modulo (*) if it has a component of degree $\geq\ n+1$, and let $D$ be a
diagram whose largest component $C$ has degree $n$.  (Of course, $n \geq 2$, and by Lemmas
\ref{L:base1} and \ref{L:base2} we can assume $k > n+2$.)  Without loss of generality,
$C$ has endpoints colored $1,2,...,n+1$.
The inner induction, which is the rest of the proof, is on $\sum_{a=n+2}^k{m(D;1,a)}$, inducting
among diagrams with a component with endpoints colored $1,2,...,n+1$.
Without loss of generality, as before, we can assume that $C$ has a branch as shown:
$$C:\ \ \begin{matrix} \bar{C} \\ | \\ | \\ 1-----2 \end{matrix}$$
We apply (*) using the color 1 and find, after removing boring diagrams and those which are trivial
by the first inductive hypothesis, that $D + m(D;1,2)D + \sum_{a=n+2}^k{m(D;1,a)D_a} = 0$, where $D_a$ is
the same as $D$ except that:
\begin{itemize}
    \item  $C$ has been replaced by a component $C_a$ identical to it except that the
endpoint colored 2 is now colored $a$ (so $\bar{C_a} = \bar{C}$)
    \item  A line segment with endpoints colored 1 and $a$ has been replaced by a line segment
with endpoints colored 1 and 2.  In other words, $m(D_a;1,a) = m(D;1,a)-1$ and $m(D_a;1,2) =
m(D;1,2)+1$.
\end{itemize}
We will denote this as shown below:
$$D_a:\ \ \begin{matrix} \bar{C} \\ | \\ | \\ 1-----a \end{matrix} \  (1,a) \rightarrow (1,2)$$
Notice that if $\sum_{a=n+2}^k{m(D;1,a)} = 0$, then $m(D;1,a) = 0$ for each $a$, since these are
all non-negative integers.  In this case, $D + m(D;1,2)D = 0$, and hence $D = 0$.  This
proves the base case of the second induction.
As in Lemma~\ref{L:base2}, we use the IHX relation to decompose $D_a = \sum_{i=3}^{n+1}{\pm D_a^i}$,
where the analogue $C_a^i$ of $C_a$ in $D_a^i$ has a branch as shown, and the other components of the
diagram are the same as $D_a$:
$$D_a^i:\ \ \begin{matrix} \bar{C_a^i} \\ | \\ | \\ i-----a \end{matrix} \ (1,a) \rightarrow (1,2)$$
Note that, aside from having endpoints of the same colors, $C_a^i$ looks nothing like $C_a$.
$\bar{C_a^i}$ has endpoints colored $1,3,4,...,i-1,i+1,...,n+1$.  We will keep this in mind.
Now we apply (*) to $D_a^i$, using the color $i$.  In the pictures we use to describe the various
diagrams that we produce in what follows, we will just be showing how the diagrams differ from
$D_a^i$.  This will involve showing how $C_a^i$ has been altered, and which line segments have
been added or removed.  At each stage, we will eliminate without comment those diagrams which
are either boring or are trivial modulo (*) by our first inductive hypothesis (i.e. have components
of degree greater than $n$).  We obtain the relation:
$$D_a^i + m(D_a^i;i,a)D_a^i + m(D_a^i;i,2)D_{a2}^i +
\sum_{\substack{n+2\leq b\leq k \\ b\neq a}}{m(D_a^i;i,b)D_{ab}^i} = 0$$
where:
$$D_{a2}^i:\ \ \begin{matrix} \bar{C_a^i} \\ | \\ | \\ i-----2 \end{matrix} \ (i,2) \rightarrow (i,a)$$
$$D_{ab}^i:\ \ \begin{matrix} \bar{C_a^i} \\ | \\ | \\ i-----b \end{matrix} \ (i,b) \rightarrow (i,a)$$
Notice that $D_{a2}^i$ has a component of degree $n$ with endpoints colored $1,2,...,n+1$, and
$m(D_{a2}^i;1,a) = m(D_a;1,a) = m(D;1,a)-1$ (since $i\neq 1$), so $D_{a2}^i$ is trivial by the
second inductive hypothesis.  So we can rewrite the relation as:
$$(1 + m(D_a^i;i,a)D_a^i) + \sum_{\substack{n+2\leq b\leq k \\ b\neq a}}{m(D_a^i;i,b)D_{ab}^i} = 0$$
Next we apply (*) to $D_{ab}^i$, using the color $b$, and find that:
$$D_{ab}^i + m(D_{ab}^i;i,b)D_{ab}^i + m(D_{ab}^i;2,b)D_{ab2}^i +
\sum_{\substack{n+2\leq c\leq k \\ c\neq b}}{m(D_{ab}^i;b,c)D_{abc}^i} = 0$$
where:
$$D_{ab2}^i:\ \ \begin{matrix} \bar{C_a^i} \\ | \\ | \\ 2-----b \end{matrix} \
\begin{matrix} (i,b) \rightarrow (i,a) \\ (2,b) \rightarrow (i,b) \end{matrix}
\Rightarrow (2,b) \rightarrow (i,a)$$
$$D_{abc}^i:\ \ \begin{matrix} \bar{C_a^i} \\ | \\ | \\ c-----b \end{matrix} \
\begin{matrix} (i,b) \rightarrow (i,a) \\ (c,b) \rightarrow (i,b) \end{matrix}
\Rightarrow (c,b) \rightarrow (i,a)$$
Now we apply (*) to $D_{ab2}^i$ using the color 2, and to $D_{abc}^i$, using the color $c$.  We
get two relations:
$$D_{ab2}^i + m(D_{ab2}^i;2,b)D_{ab2}^i + m(D_{ab2}^i;2,i)D_{ab2i}^i +
\sum_{\substack{n+2\leq c\leq k \\ c\neq b}}{m(D_{ab2}^i;2,c)D_{ab2c}^i} = 0$$
$$D_{abc}^i + m(D_{abc}^i;b,c)D_{abc}^i + m(D_{abc}^i;2,c)D_{abc2}^i + m(D_{abc}^i;i,c)D_{abci}^i$$
$$+ \sum_{\substack{n+2\leq d\leq k \\ d\neq b,c}}{m(D_{abc}^i;c,d)D_{abcd}^i} = 0$$
where:
$$D_{ab2i}^i:\ \ \begin{matrix} \bar{C_a^i} \\ | \\ | \\ 2-----i \end{matrix} \
\begin{matrix} (2,b) \rightarrow (i,a) \\ (2,i) \rightarrow (2,b) \end{matrix}
\Rightarrow (2,i) \rightarrow (i,a)$$
$$D_{ab2c}^i:\ \ \begin{matrix} \bar{C_a^i} \\ | \\ | \\ 2-----c \end{matrix} \
\begin{matrix} (2,b) \rightarrow (i,a) \\ (2,c) \rightarrow (2,b) \end{matrix}
\Rightarrow (2,c) \rightarrow (i,a)$$
$$D_{abc2}^i:\ \ \begin{matrix} \bar{C_a^i} \\ | \\ | \\ c-----2 \end{matrix} \
\begin{matrix} (c,b) \rightarrow (i,a) \\ (2,c) \rightarrow (c,b) \end{matrix}
\Rightarrow (2,c) \rightarrow (i,a)$$
$$D_{abci}^i:\ \ \begin{matrix} \bar{C_a^i} \\ | \\ | \\ c-----i \end{matrix} \
\begin{matrix} (c,b) \rightarrow (i,a) \\ (i,c) \rightarrow (c,b) \end{matrix}
\Rightarrow (i,c) \rightarrow (i,a)$$
$$D_{abcd}^i:\ \ \begin{matrix} \bar{C_a^i} \\ | \\ | \\ c-----d \end{matrix} \
\begin{matrix} (c,b) \rightarrow (i,a) \\ (c,d) \rightarrow (c,b) \end{matrix}
\Rightarrow (c,d) \rightarrow (i,a)$$
We make several observations, by antisymmetry:
\begin{itemize}
    \item  $D_{ab2i}^i = -D_{a2}^i = 0$.
    \item  $D_{ab2c}^i = -D_{abc2}^i$.
    \item  $D_{abci}^i = -D_{ac}^i$.
    \item  $D_{abcd}^i = D_{adc}^i = -D_{acd}^i$.
\end{itemize}
Now that we have these recursive relations, we can plug them into our various equations.  We will
use the following equalities:
\begin{eqnarray*}
m(D_{ab}^i;i,b) + 1 & = & m(D_a^i;i,b) \\
m(D_{ab2}^i;2,b) + 1 & = & m(D_{ab}^i;2,b) = m(D_a^i;2,b) \\
m(D_{abc}^i;b,c) + 1 & = & m(D_a^i;b,c) \\
\end{eqnarray*}
And for all the other coefficients we have:
$$m(D_*^i;x,y) = m(D_a^i;x,y)$$
For convenience, we will write $m(x,y) = m(D_a^i;x,y)$ in what follows:
\begin{eqnarray*}
(m(i,a)+1)D_a^i & = & \sum_{\substack{n+2\leq b\leq k \\ b\neq a}}{-m(i,b)D_{ab}^i} \\
 & = & \sum_{\substack{n+2\leq b\leq k \\ b\neq a}}{-(m(D_{ab}^i;i,b)+1)D_{ab}^i} \\
 & = & \sum_{\substack{n+2\leq b\leq k \\ b\neq a}}{\left({m(2,b)D_{ab2}^i +
    \sum_{\substack{n+2\leq c\leq k \\ c\neq b}}{m(b,c)D_{abc}^i}}\right)} \\
\end{eqnarray*}
Note that:
\begin{eqnarray*}
m(2,b)D_{ab2}^i & = & (m(D_{ab2}^i;2,b)+1)D_{ab2}^i \\
 & = & \sum_{\substack{n+2\leq c\leq k \\ c\neq b}}{-m(2,c)D_{ab2c}^i} \\
m(b,c)D_{abc}^i & = & (m(D_{abc}^i;b,c)+1)D_{abc}^i \\
 & = & -m(2,c)D_{abc2}^i - m(i,c)D_{abci}^i - \sum_{\substack{n+2\leq d\leq k \\ d\neq c,b}}
    {m(c,d)D_{abcd}^i} \\
\end{eqnarray*}
Therefore:
$$m(2,b)D_{ab2}^i + \sum_{\substack{n+2\leq c\leq k \\ c\neq b}}{m(b,c)D_{abc}^i} = $$
$$\sum_{\substack{n+2\leq c\leq k \\ c\neq b}}{\left({-m(2,c)(D_{ab2c}^i + D_{abc2}^i)
    - m(i,c)D_{abci}^i - \sum_{\substack{n+2\leq d\leq k \\ d\neq c,b}}{m(c,d)D_{abcd}^i}}\right)} = $$
$$\sum_{\substack{n+2\leq c\leq k\\ c\neq b}}{\left({m(i,c)D_{ac}^i +
    \sum_{\substack{n+2\leq d\leq k \\ d\neq c,b}}{m(c,d)D_{acd}^i}}\right)}$$
We plug this back in above to find:
$$(m(i,a)+1)D_a^i = \sum_{\substack{n+2\leq b\leq k \\ b\neq a}\ }
    {\sum_{\substack{n+2\leq c\leq k \\ c\neq b}}{\left({m(i,c)D_{ac}^i +
    \sum_{\substack{n+2\leq d\leq k \\ d\neq c,b}}{m(c,d)D_{acd}^i}}\right)}}$$
We notice that:
\begin{eqnarray*}
\sum_{\substack{n+2\leq c\leq k \\ c\neq b}\ }{\sum_{\substack{n+2\leq d\leq k \\ d\neq c,b}}
    {m(c,d)D_{acd}^i}} & = & \frac{1}{2}\sum_{\substack{n+2\leq c,d\leq k \\ c\neq d
    \\ c,d\neq b}}{m(c,d)(D_{acd}^i + D_{adc}^i)} \\
 & = & \frac{1}{2}\sum_{\substack{n+2\leq c,d\leq k \\ c\neq d \\ c,d\neq b}}{m(c,d)
    (D_{acd}^i - D_{acd}^i)} \\
 & = & 0 \\
\end{eqnarray*}
Therefore:
$$(m(i,a)+1)D_a^i = \sum_{\substack{n+2\leq b\leq k \\ b\neq a}\ }
    {\sum_{\substack{n+2\leq c\leq k \\ c\neq b}}{m(i,c)D_{ac}^i}}$$
Returning to our first equation above, we have (simply replacing $b$ by $c$ in the second
equality):
\begin{eqnarray*}
(m(i,a)+1)D_a^i & = & \sum_{\substack{n+2\leq b\leq k \\ b\neq a}}{-m(i,b)D_{ab}^i} \\
 & = & \sum_{\substack{n+2\leq c\leq k \\ c\neq a}}{-m(i,c)D_{ac}^i} \\
 & = & \left({\sum_{\substack{n+2\leq c\leq k \\ c\neq b}}{-m(i,c)D_{ac}^i}}\right)
    + m(i,a)D_{aa}^i - m(i,b)D_{ab}^i \\
\end{eqnarray*}
Since $D_{aa}^i = D_a^i$, we can cancel and rearrange terms to write:
$$\sum_{\substack{n+2\leq c\leq k \\ c\neq b}}{m(i,c)D_{ac}^i} = -D_a^i - m(i,b)D_{ab}^i$$
Therefore:
\begin{eqnarray*}
(m(i,a)+1)D_a^i & = & \sum_{\substack{n+2\leq b\leq k \\ b\neq a}}{-D_a^i - m(i,b)D_{ab}^i} \\
 & = & -(k-n-2)D_a^i + (m(i,a)+1)D_a^i \\
\end{eqnarray*}
$$\Rightarrow (k-n-2)D_a^i = 0$$
Since $n < k-2$ (the cases when $n=k-1,k-2$ were dealt with in the lemmas), $k-n-2 \neq 0$; so we can
conclude that $D_a^i$ is trivial modulo (*).  Hence, $D_a$ and, ultimately, $D$ are also trivial modulo
(*).  This completes the induction and the proof.  $\Box$
This theorem tells us that the only elements of $B^{hsl}$ which are {\it not} in the kernel of the
relation (*) are unitrivalent diagrams all of whose components are of degree 1 (i.e. line segments).
Restricted to the space generated by these elements, (*) is clearly trivial, so $B^{hl}$ is in fact
simply the polynomial algebra over the reals generated by these unitrivalent diagrams (since (*)
is trivial on this space, $B^{hl}$ inherits a multiplication from $B^{hsl}$).  We formalize this
as a corollary:
\begin{cor} \label{C:algebra}
$B^{hl}(k)$ (and hence $A^{hl}(k)$) is isomorphic to the algebra
${\mathbb R}[x_{ij}]$,
where each $x_{ij}$ is of degree 1, and $1\leq i < j\leq k$.
\end{cor}
It is well-known that the only finite-type link homotopy invariants of degree 1 are the pairwise
linking numbers of the components, so we conclude:
\begin{cor} \label{C:reduction}
The pairwise linking numbers of the components of a link are the only real-valued
finite type link homotopy invariants of the link.
\end{cor}
{\bf Remark:}  Bar-Natan~\cite{bn2} has shown that the Milnor $\mu$-invariants are finite
type homotopy invariants for {\it string} links.  However, the analogous
$\bar{\mu}$-invariants for links have indeterminacies arising from the fact that many string
links can close up to give the same link (up to homotopy).  As a result, these invariants are
only well-defined modulo the values of lower-order $\mu$-invariants.  This keeps us from
being able to extend the invariants to singular links, since two links which differ by a
crossing change may have entirely different lower-order invariants, and so their
$\bar{\mu}$-invariants may have values lying in completely different groups.  So there
is no way to interpret these invariants as finite type invariants in the usual way.

\section{Acknowledgements}
I wish to particularly thank Robert Schneiderman, who asked me the question which motivated this
paper.  I also wish to thank my advisor, Robion Kirby, for his advice and support, and all the
members of the Informal Topology Seminar at Berkeley.  Finally, I want to thank the
anonymous referee for his or her many valuable comments.

\end{document}